\documentclass[11pt]{article}

\usepackage{amssymb, amsmath, latexsym}

\newtheorem{lem}{Lemma}[section]
\newtheorem{thm}{Theorem}[section]
\newtheorem{cor}{Corollary}[section]

\newenvironment{proof}[1]{ \textbf{#1}}{\hfill $\Box$\vspace{2ex}}
\newenvironment{remark}[1]{\vspace{2ex}\hspace{-\parindent}\textbf{Remark {#1}.}}{\hfill $\Box$\vspace{2ex}}
\newenvironment{observation}[1]{\vspace{2ex}\hspace{-\parindent}\textbf{{#1}}}{\hfill $\Box$\vspace{2ex}}

\title{Some families of special Lagrangian tori}
\author{Diego Matessi \thanks
{Mathematics Institute, University of Warwick CV4 7AL, Coventry, UK.
e-mail: diego@maths.warwick.ac.uk}}

\begin{document}

\maketitle

\newcommand{\dery}[2]{\ensuremath{\frac{\partial {#1}}{\partial y_{#2}}}}
\newcommand{\derx}[2]{\ensuremath{\frac{\partial {#1}}{\partial x_{#2}}}}
\newcommand{\delx}[1]{\ensuremath{\frac{\partial}{\partial x_{#1}}}}
\newcommand{\dely}[1]{\ensuremath{\frac{\partial}{\partial y_{#1}}}}
\newcommand{\mbar}{\ensuremath{\overline{M}}}
\newcommand{\uc}{\ensuremath{U_{\complex{}}}}
\newcommand{\reals}[1]{\ensuremath{\mathbb{R}^{#1}}}
\newcommand{\integ}[1]{\ensuremath{\mathbb{Z}^{#1}}}
\newcommand{\nat}{\ensuremath{\mathbb{N}}}
\newcommand{\complex}[1]{\ensuremath{\mathbb{C}^{#1}}}
\newcommand{\projc}[1]{\ensuremath{\mathbb{C} \mathbb{P}^{#1}}}
\newcommand{\quaternions}[1]{\ensuremath{\mathbb{H}^{#1}}}
\newcommand{\cayley}[1]{\ensuremath{\mathbb{O}^{#1}}}
\newcommand{\haus}[1]{\ensuremath{\mathcal{H}^{#1}}}
\newcommand{\Img}{\operatorname{Im}}
\newcommand{\Rl}{\operatorname{Re}}
\newcommand{\Vl}{\operatorname{Vol}}
\newcommand{\Rc}{\operatorname{Ric}}
\newcommand{\Hm}{\operatorname{Hom}}
\newcommand{\spt}{\operatorname{spt}}
\newcommand{\Hes}{\operatorname{Hess}}
\newcommand{\Identity}{\operatorname{Id}}
\newcommand{\Ed}{\operatorname{End}}
\newcommand{\tr}{\operatorname{Tr}}
\newcommand{\inner}[2]{ \langle {#1}, {#2} \rangle}
\newcommand{\Lapl}{\Delta^{\nabla}}

\begin{abstract}
We give a simple proof of the local version of Bryant's
re\-sult~\cite{bry:calembed}, sta\-ting that any 3-dimensional 
Riemannian manifold can be isome\-tri\-cal\-ly embedded as a 
special Lagrangian submanifold in a Calabi-Yau manifold.
We then refine the theorem proving that a certain class of one-parameter
families of metrics on a 3-torus can be isometrically embedded in a
Calabi-Yau manifold as a one-parameter family of 
special Lagrangian submanifolds.
Two applications of our results show how the geometry of moduli space of
3-dimesional special Lagrangian submanifolds  differs considerably 
from the 2-dimensional one.
First of all, applying our first theorem and a construction due to Calabi 
we show that nearby elements of the local moduli space of a special 
Lagrangian 3-torus can intersect themselves. 
Secondly, we use our examples of one-parameter families to show that the 
semi-flat metric on the mirror manifold proposed by Hitchin in 
\cite{hitch:msslag} is not necessarily Ricci-flat in dimension 3.
\end{abstract}

\section{Introduction}
Many interesting speculations have been made about the role 
special Lagrangian submanifolds should play in understanding
the geometry of Calabi-Yau manifolds and of Mirror Symmetry.
Unfortunately the lack of examples has allowed few of these
to be proved.
Only recently has the number of new constructions finally begun to increase.
For years, in fact, the only examples known were the ones
appearing in the foundational paper by Harvey and Lawson~\cite{hl:cal},  where special
Lagrangian submanifolds were defined for the first time.
Our paper participates in the quest for examples.
We propose a new way to construct special Lagrangian
submanifolds and one-parameter families of these and we relate them
to some of the speculations which have been made about them. 
Let's first recall some definitions. For us, a Calabi~-~Yau manifold will be a triple
$(\mbar, \Omega, \omega)$ where $\mbar$ is a complex  n-dimensional
manifold, $\Omega$ a nowhere-vanishing holomorphic n-form on $\mbar$ and $\omega$
a K\"{a}hler form related to $\Omega$ by
\begin{equation} \label{cy:str}
 \omega^{n} = c \Omega \wedge \overline{\Omega}, 
\end{equation}
for some constant c. By Yau's proof of the Calabi conjecture, this triple
can be constructed on any compact K\"{a}hler manifold with trivial canonical bundle. 
The K\"{a}hler metric $\omega$ is Ricci-flat. An n-dimensional submanifold $M$ is called
special Lagrangian (sometimes abbreviated sLag) if it satisfies:
\[ \Rl{\Omega}_{|M} = \Vl_{M}, \]
where $\Vl_{M}$ denotes the volume form on $M$. Equivalently, $M$ is special
Lagrangian if and only if it satisfies the following:
\begin{eqnarray}
\Img{\Omega}_{|M} &= & 0,  \label{slag:def} \\ 
\omega_{|M} &= & 0.       \label{slag:def1}
\end{eqnarray}
In\, this\, pa\-per\, we\, will\, ve\-ry\, often\, re\-fer\, to\, the   
work of three au\-thors:
McLean~\cite{mclean:deform}, Hitchin~\cite{hitch:msslag} and 
Gross~\cite{mgross:slfibtop,mgross:slfibgeom,mgross:topms}.
We briefly describe here their results. 
Given a special Lagrangian submanifold $M$, McLean proved
that the moduli space of nearby special Lagrangian submanifolds
can be identified with a smooth submanifold $\mathcal{M}$ of $\Gamma(\nu(M))$,
the space of sections of $\nu(M)$, the normal bundle of $M$. The dimension of
$\mathcal{M}$ is $b_{1}(M)$, the first Betti number of $M$. In fact, through the
map $V \rightarrow (JV)^{\flat}$ (cfr. end of section for notation),
which identifies a section $V$ in $\Gamma(\nu(M))$ with a section
in $\Omega^{1}(M)$, $\mathcal{M}$ can be viewed inside $\Omega^{1}(M)$
and its tangent space at $M$ turns out to be  the vector space of harmonic one-forms
on $M$. 
In practice, the latter means that if we take a variation of $M$
through special Lagrangian 
submanifolds with variational vector field $V$, then $(JV)^{\flat}$
is a harmonic one-form. In particular, if  $M$ is a torus with non-vanishing harmonic one-forms,
then McLean's result implies that a whole open set of $\mbar$
around $M$ is fibred by special Lagrangian tori. On $\mathcal{M}$ there is also
a natural metric which is the standard $L^{2}$ norm of one-forms.

In \cite{syz:slfib} the three authors conjectured,
in what is now called the SYZ-conjecture, a geometric construction of Mirror Symmetry.
Here, on purely physical grounds, they argued that if $\mbar$ 
is near some boundary point of its complex moduli
space then it should be possible to fibre it through special Lagrangian 
tori, some of which may be singular.
The mirror manifold of $\mbar$, in the sense of Mirror Symmetry, 
is obtained by dualizing this fibration.  
Some mathematical aspects the conjecture were described and
investigated by Hitchin~\cite{hitch:msslag} and 
Gross~\cite{mgross:slfibgeom, mgross:slfibtop, mgross:topms}.
First Hitchin showed how $\mathcal{M}$ can be naturally identified
with an open subset of $H^{1}(M,\reals{})$ or of $H^{n-1}(M,\reals{})$
and explained how the two identifications are dual to each other. 
According to the SYZ-conjecture, in the case $M$ is a torus,
a local candidate for the mirror of $\mbar$ is the space
\[ \mathcal{X} = \mathcal{M} \times H^{1}(M, \reals{}/ \integ{}). \]
This is a torus fibration over $\mathcal{M}$.
The problem is to find, possibly in a natural way, a
Calabi-Yau structure on this fibration such that the fibres are
special Lagrangian tori. Using the identifications above, Hitchin 
explained how to construct an integrable
complex structure, a K\"{a}hler form and a holomorphic n-form on $\mathcal{X}$.
This metric is often called the
semi-flat metric. He then showed that these forms give a Calabi~-~Yau structure,
i.e. they are related by (\ref{cy:str}),
if and only if $\mathcal{M}$ satisfies a certain condition.
While this condition is known to be satisfied in the 2-dimensional case
(see for example Hitchin~\cite{hitch:msclag}),
it is one of the results of this paper that in general it is not in dimension 3. 

Gross dealt with the more global aspects of the SYZ construction by
treating the problem of how to include singular special Lagrangian 
fibres in the above picture. In fact, on the basis of the
topological consequences of Mirror Symmetry, he gave a
conjectural description of the singular fibres which are
expected to appear and explained how to dualize them. 
This construction is completely understood for K3 surfaces,
where special Lagrangian fibrations are just elliptic fibrations
with a different complex structure.

Parallel to these speculative aspects of special Lagrangian
geometry, there has been the attempt to produce examples. 
After the Harvey and Lawson ones, Bryant~\cite{bry:realslice} and Kobayashi~\cite{kob:realslice}
showed how to construct special Lagrangian
tori as totally real submanifolds of subvarieties of $\projc{n}$. 
Lately, many examples of special Lagrangian fibrations where constructed on
complete Calabi-Yau manifolds by 
Goldstein~\cite{goldstein:calfib, goldstein:cftorus, goldstein:sltorusact}.
In~\cite{mgross:slexmp} Gross used similar ideas to Goldstein's to construct
special Lagrangian fibrations on $\complex{n}/G$, where $G$ is a finite
abelian subgroup of $SU(n)$. 
More recently Haskins~\cite{haskins:slcones} found more special
Lagrangian cones in in $\complex{3}$. His construction was subsequently
generalized by Joyce~\cite{joyce:slsymmetries, joyce:slquadrics}, who
also provided other examples which are not cones. 

The results of this paper overlap in part 
with those obtained by Bryant~\cite{bry:calembed}. He proved 
that any real-analytic, 3-dimensional Riemannian manifold
$(M,g)$ with real-analytic metric $g$ can be isometrically embedded in some
Calabi-Yau manifold $\mbar$. His proof used Cartan-K\"{a}hler theory,
which requires the problem to be translated into one of existence of
integral submanifolds of a differential ideal.  
Our first result (Theorem~\ref{splag:embedd}) is the local version 
of the same theorem, but the proof is simpler
and is global in the case of the torus. We prove the following: given any
pair $(U,g)$ where $U$ is some open set in $\reals{3}$ and $g$ a metric,
we can isometrically embed $U$ as a special Lagrangian submanifold of
some Calabi-Yau manifold $\mbar$. 
Our proof, as well as being simple, has other advantages. First of all we show
that the complex structure of $\mbar$ around $U$ is in some sense unique and
can be dealt with very concretely with a suitable choice of coordinates. 
Hence we prove that also the holomorphic n-form is unique, in fact it is literally
the holomorphic extension of the volume form on $U$.
Finally, we write the equations for the Ricci-flat K\"{a}hler metric 
and show that a solution always exist with three successive applications
of the Cauchy-Kowalesky theorem.
Using this result and the construction by Calabi of metrics on
the 3-torus which admit harmonic one-forms with zeroes we show
that there are special Lagrangian 3-tori which can intersect
elements of the moduli space of its deformations. This did not
happen in dimension 2.

The structure of the proof of our first result leads to an immediate 
refinement.  In fact we show (Theorem~\ref{onepar:constr}) that if a
one-parameter family of metrics on a 3-torus satisfies certain simple conditions, 
then it can always be realized as a one-parameter family of special Lagrangian tori
in a Calabi-Yau manifold. The set of one-parameter families thus constructed
is quite rich and provides us with many examples. 
Some of these also show that the condition required for Hitchin's metric to yield a
Calabi-Yau structure is in fact not satisfied. This leads to the question
of how can one find such a structure.

\textbf{Notations.} When working in $\complex{n}$ 
complex coordinates are always denoted by $(z_{1}, \ldots, z_{n})$,
and real coordinates by $(x_{1},\ldots,x_{n}, y_{1}, \ldots, y_{n})$, where
$z_{k} = x_{k}+iy_{k}$. Sometimes we will use $x$ (or $y$) as short
for $(x_{1}, \ldots, x_{n})$ (or $(y_{1}, \ldots, y_{n})$).
The letter $J$ is always used to denote the almost complex structure.
The superscript $(V)^{\flat}$ stands for the
element in $T^{\ast}M$ cor\-res\-pon\-ding to $V$ under the identification
of $TM$ and $T^{\ast}M$ induced by the metric.
As usual $\star: \Omega^{k}(M) \rightarrow \Omega^{n-k}(M)$ denotes
the Hodge-star operator between forms.
We follow the convention that given the coefficients of an invertible
matrix $g_{ij}$, the terms $g^{ij}$ denote the coefficient of the 
inverse matrix. 

\textbf{Acknowledgements.} The author wishes to thank his thesis 
advisors Mark Gross and Mario J. Micallef for the invaluable help
he received from them and for introducing him to this fascinating
subject. He also wishes to thank Rita Gaio and Luca Sbano for some 
very usefull discussions.

\section{Complexifications}
Given a real-analytic, n-dimensional manifold $M$, a 
\textbf{complexification} of $M$ is an n-dimensional complex manifold
$\overline{M}$ together with a real analytic embedding 
$\iota : M \rightarrow \overline{M}$ such that for every
$p \in \overline{M}$ there exist holomorphic coordinates 
$(z_{1}, \ldots, z_{n})$ on a neighborhood $U$ of $p$ such that 
$q \in U \cap \iota(M)$ if and only if $\Img (z_{i}(q)) = 0, \ \ i=1, \ldots, n$. 

\begin{observation}{Example 1.}
Given an open set $U \subseteq \reals{n}$, identify it with a 
subset of $\complex{n}$ through the standard inclusion of
$\reals{n}$ in $\complex{n}$ as the real part. An open 
neighborhood $\uc$ of $U$ such that $\Rl({\uc})=U$ will
be called a \textbf{standard complexification} of $U$.
So, $\mbar$ being a complexification of $M$ means that,
locally, the pair $(\mbar, M)$ is holomorphic to the
pair $(\uc,U)$.
\end{observation}

\begin{observation}{Example 2.}
Let $M$ be the standard n-torus $\reals{n}/ \integ{n}$
and $\iota$ its obvious inclusion in $\complex{n}/ \integ{n}$,
where $\integ{n}$ acts through translations on the real part. 
Then $(\complex{n}/ \integ{n}, \iota)$ is a complexification
of $M$. It will be referred to as a \textbf{standard complexification}
of the n-torus.
\end{observation}
 
Bruhat and Whitney~\cite{whbr:cx} proved the following:
\begin{thm} \textbf{(Bruhat, Whitney)}
Any paracompact, real-analytic manifold $M$ admits a complexification.
Moreover if $(\overline{M}_{1},\iota_{1})$ and $(\overline{M}_{2}, \iota_{2})$ are
two complexifications of $M$, then there exist neighborhoods $V_{i}$ of 
$\iota_{i}(M)$, $i=1,2$, and a biholomorphism $F: V_{1} \rightarrow V_{2}$ extending 
$\iota_{2} \circ \iota_{1}^{-1}$.
\label{bruhat:whitney}
\end{thm}
They also showed that there exists an antiholomorphic 
involution $\sigma: \overline{M} \longrightarrow \overline{M}$ which has $M$ as
the set of its fixed points.
Identify $\iota(M)$ with $M$. We say that $M$ is a \textbf{totally real} submanifold
of a complex manifold $\overline{M}$ if $J(T_{p}M)$ is transversal to $T_{p}M$, for
for every $p \in M$, where $J$ is the complex structure on $\overline{M}$. If 
$\overline {M}$ is a complexification of $M$ then $M$ is obviously a totally real
submanifold of $\overline{M}$. The converse is also true:
\begin{lem}
Let $\iota : M \rightarrow \overline{M}$ be a real-analytic embedding 
of $M$ as a totally real submanifold of the complex manifold $\overline{M}$.
Then $(\overline{M}, \iota)$ is
a complexification of M. 
\label{cx:totrl}
\end{lem}
\begin{proof}{Proof.} Let $p \in M$. We can assume w.l.o.g.
$\overline{M} = \complex{n}$, $p=0$ and 
$T_{p}M = \{ \Img (z_{i}) = 0, \ i=1, \ldots , n \}$. 
Then there exists a neighborhood $V \subset \complex{n}$ of $0$ 
and a real-analytic map $f: \Rl(V) \rightarrow \reals{n}$ such that 
$V \cap M = \{ x + if(x), x \in V \}$. Extend $f$ to a holomorphic function
$\tilde{f} : \tilde{V} \rightarrow \complex{n}$,
where $\tilde{V}$ is some neighborhood of $\Rl(V)$ in $\complex{n}$. Define
$\tilde{F} : \tilde{V} \rightarrow \complex{n}$ by $\tilde{F}(z) = z + i \tilde{f}(z)$,
then $\tilde{F}$ is a biholomorphism near $0$ and $F= \tilde{F}^{-1}$ gives
the complex coordinates with the required property.
\end{proof}

In particular we have the following:
\begin{cor}
Let $\overline{M}_{1}$ be a
K\"{a}hler manifold and  and $\iota_{1}:M \rightarrow \overline{M}_{1}$
a real-analytic embedding of $M$ as a Lagrangian submanifold. If $(\overline{M}_{2},
\iota_{2})$ is a complexification of $M$, then there exist neighborhoods $V_{i}$ 
of $\iota_{i}(M)$ and a biholomorphism $F: V_{1} \rightarrow V_{2}$ extending 
$\iota_{2} \circ \iota_{1}^{-1}$.
\label{cx:lag}
\end{cor}
\begin{proof}{Proof.} It follows immediately from Theorem~\ref{bruhat:whitney}
and Lemma~\ref{cx:totrl} since Lagrangian submanifolds are totally real.
\end{proof}

Notice that, since special Lagrangian submanifolds
are minimal, they are also real-analytic. Hence Corollary~\ref{cx:lag}
applies when $M$ is a special Lagrangian submanifold. In particular
if $\phi: U \rightarrow M$ is a real-analytic coordinate
chart, it can be extended to a holomorphic chart 
$\phi_{\complex{}}: \uc \rightarrow \mbar$. Also, in the
case $M$ is the n-torus and $(\complex{n}/ \integ{n}, \iota)$
its standard complexification, then any special Lagrangian embedding 
$\tau: M \rightarrow \mbar$ can be extended to
a holomorphic chart $F:\uc \rightarrow \mbar$, where
$\uc$ is a sufficiently small neighborhood of $M$ in $\complex{n}/ \integ{n}$.

\section{Local isometric special Lagrangian embeddings} \label{isom:slag}
Now let $(U,g)$ be an open neighborhood of $0 \in \reals{3}$ together with
a Riemannian metric $g=(g_{ij})$. We look for isometric embeddings
of $(U,g)$ as a special Lagrangian submanifold of some Calabi-Yau $\mbar$. 
From the results in the previous section
we may assume w.l.o.g. that $\mbar= \uc$ for some standard complexification $\uc$.
Remember that $\uc$ is a subset of $\complex{n}$, so we can use the standard
complex coordinates $(z_{1}, \ldots , z_{n})$.
 We will prove the following:
\begin{thm} \label{splag:embedd}
On some standard complexification $\uc$ of $U$ we can find a unique holomorphic
n-form $\Omega$ and at least one K\"{a}hler form $\omega$ satisfying the
following properties:
\begin{enumerate}
\item $\omega^{3}/3! = -(i/2)^{3} \Omega \wedge \overline{\Omega}$, 
\item the induced metric on $U$ is $g$,
\item $\Omega_{|U} = \Vl_{U}$.
\end{enumerate}
\end{thm}
The first condition is just equation (\ref{cy:str}) from the Introduction, 
with a choice of the constant $c$. Conditions 2 and 3 make 
$(U,g)$ isometrically embedded
in $(\uc, \Omega, \omega)$ as a special Lagrangian submanifold. In what follows
we will denote by $h=(h_{ij})$ the hermitian metric associated with $\omega$. 
Part of the theorem is proved by the next lemma:
\begin{lem}
There exists a unique $\Omega$ on $\uc$ satisfying conditions (1)-(3) above. In fact,
in standar coordinates, $\Omega$ must be
\[ \Omega = \Gamma_{g} (z) dz_{1} \wedge dz_{2} \wedge dz_{3}, \]
where $\Gamma_{g}$ denotes the holomorphic extension of 
$\sqrt{g} = \sqrt{\det (g_{ij})}$, the coefficient of $\Vl_{U}$.
\end{lem}
\begin{proof}{Proof.}
Certainly we can write
\[ \Omega = f(z) dz_{1} \wedge dz_{2} \wedge dz_{3}, \]
for some holomorphic $f$. Let $f = \alpha + i \beta $, then condition (1) gives:
\[ \det ({h}_{ij}) = \alpha^{2} + \beta^{2}. \]
From condition (2) it follows that, along $U$, we have  $h_{ij} (x, 0) = g_{ij}(x)$, 
giving that $\det (h_{ij})(x,0) = g(x)$. Condition (3) implies that
\[ \Omega_{|U} = \alpha dx_{1} \wedge dx_{2} \wedge dx_{3} = \sqrt{g} dx_{1} \wedge dx_{2} \wedge dx_{3}. \]
Therefore we obtain that $\beta (x, 0) = 0$ and $ f(x, 0) = \alpha (x,0) = \sqrt{g}(x)$.
The only holomorphic function satisfying this is precisely $\Gamma_{g}$.
\end{proof}

\begin{proof}{Proof of Theorem~\ref{splag:embedd}}
We write the hermitian metric $h$ that we are looking for as $h=A+iB$, 
where $A=(\alpha_{ij})$ and $B=(\beta_{ij})$ are real valued matrices, symmetric 
and antisymmetric respectively. In the basis $(\delx{1},\ldots,\delx{n},\dely{1}, \ldots, \dely{n})$ 
for $T\uc$ the corresponding K\"{a}hler form can be written as a $2n \times 2n$ matrix
\[ \omega = \left( \begin{array}{cc}
                                  -B & A \\
                                  -A & -B 
                               \end{array}  \right). 
                              \]
In order to prove the theorem we need to solve the following ``initial value''
PDE problem:
\[ \left\{ \begin{array}{lr}
                   \det(h)= |\Gamma_{g}|^2 & (D)\\
                   d\omega = 0 & (C)\\
                   A(x,0) = g(x) \ \text{and} \ B(x,0) = 0 \ \text{for all} \  x \in U. & (I)
                 \end{array} \right. \]
If we do the computations explicitly we see that (D) and (C) form the following system 
of equations in the coefficients of $\omega$:
\[ \begin{array}{lr}
       \begin{array}{lcl}
         (\alpha_{22} \alpha_{33} & - & \alpha_{23}^{2} - \beta_{23}^{2})\alpha_{11}
 - \beta_{13}^{2}  \alpha_{22} \\
\ & - & \beta_{12}^{2} \alpha_{33} - \alpha_{12}^{2}\alpha_{33} -  \alpha_{13}^{2}  \alpha_{22} 
+  2 \alpha_{12} \alpha_{23} \alpha_{13} \\
          \ &- & 2\beta_{12}\beta_{23} \alpha_{13} + 2 \alpha_{12}\beta_{23}\beta_{13} 
          + 2\beta_{12} \alpha_{23}\beta_{13} = |\Gamma_{g}|^2 
        \end{array}  & (D) \\
     \ & \ \\
     \dery{\beta_{ij}}{1} = \derx{\alpha_{1j}}{i} - \derx{\alpha_{1i}}{j} & (C1) \\
     \ & \ \\
     \dery{\beta_{ij}}{2} = \derx{\alpha_{2j}}{i} - \derx{\alpha_{2i}}{j} & (C2.1) \\
     \dery{\alpha_{1k}}{2} = \dery{\alpha_{2k}}{1} + \derx{\beta_{12}}{k} & (C2.2) \\
     \ & \ \\
    \dery{\beta_{ij}}{3} = \derx{\alpha_{3j}}{i} - \derx{\alpha_{3i}}{j} & (C3.1) \\
    \dery{\alpha_{1k}}{3} = \dery{\alpha_{3k}}{1} + \derx{\beta_{13}}{k} & (C3.2) \\
    \dery{\alpha_{2k}}{3} = \dery{\alpha_{3k}}{2} + \derx{\beta_{23}}{k} & (C3.3) \\
    \ & \ \\
    \derx{\beta_{23}}{1}-\derx{\beta_{13}}{2}+\derx{\beta_{12}}{3}=0 & (C4.1) \\
    \dery{\beta_{23}}{1}-\dery{\beta_{13}}{2}+\dery{\beta_{12}}{3}=0. & (C4.2) 
    \end{array} \]
Here the index $k$ goes from 1 to 3, while $i,j$ are such that $i <  j$.

A solution is constructed in three steps: first we find one on 
$\uc^{1}= \{ (z_{1},z_{2},z_{3}) \in \uc| y_{2}=y_{3}=0 \}$, then we extend it to 
$\uc^{2}= \{ (z_{1},z_{2},z_{3}) \in \uc| y_{3}=0 \}$ and finally to the whole $\uc$. 
Notice that for the first step we need only to look at equations (D) and (C1), 
which do not involve derivatives with respect to $y_{2}$ or $y_{3}$. 
For reasons that will become apparent later we do not assume that
$A$ is symmetric. Hence, we have four equations for twelve unknowns 
(nine from $A$ and three from $B$).
We choose arbitrarily all $\alpha_{ij}$'s on $\uc^{1}$ except 
$\alpha_{11}$, with the only requirements that they satisfy the initial conditions (I), 
they are real-analytic and they can be coefficients of a metric (e.g. $\alpha_{ij}=\alpha_{ji}$).
It is now easy to see that by differentiating (D) by $y_{1}$ 
and substituting into it equations from (C1), (D) can be written in the form
\[ \begin{array}{lr}
 \dery{\alpha_{11}}{1}  =  P(x,y_{1},\alpha_{11},\beta) \derx{\alpha_{11}}{2} + 
  Q(x,y_{1},\alpha_{11}, \beta) \derx{\alpha_{11}}{3} + R(x,y_{1}, \alpha_{11}, \beta), & (D^{\prime})
   \end{array}      \]
where $P,Q$ and $R$ are real-analytic coefficients, which depend on the 
way we arbitrarily extended the other $\alpha_{ij}$'s. Notice that this is possible
also because, with the given initial conditions, the coefficient of
$\alpha_{11}$ in (D) is different from zero near $U$.
Now equations ($D^{\prime}$) and (C1) are four equations in the four unknowns 
$\alpha_{11}, \beta_{12}, \beta_{13}, \beta_{23}$ of the type whose solution 
is guaranteed to exist uniquely (at least locally) by the Cauchy-Kowalesky theorem 
(as stated for example in Spivak \cite[Section 10.5]{spivak:dg}). The solution will also satisfy
equation (C4.1). In fact this is demonstrated by differentiating (C1), $i=1$, $j=2$ by
$x_{3}$; (C1), $i=1$, $j=3$ by $x_{2}$ and (C1), $i=2$, $j=3$ by $x_{1}$. From the 
results it follows that
\[\dely{1} \left( \derx{\beta_{23}}{1}-\derx{\beta_{13}}{2}+\derx{\beta_{12}}{3} \right) = 0 \]
on $\uc^{1}$. This shows that since equation (C4.1) holds on $U$ it holds everywhere also on $\uc^{1}$.  

The second step is similar. We now extend this solution to $\uc^{2}$ by looking at
equations (D) and the group (C2). This time we have seven equations for twelve unknowns.
We arbitrarily extend $\alpha_{33}$ and $\alpha_{23}=\alpha_{32}$ as before. 
Then, for the symmetry of $A$, we also impose $\alpha_{12}= \alpha_{21}$ and $\alpha_{13}= \alpha_{31}$.
Differentiating (D) by $y_{2}$, again we see that we can reduce the system to 
one which is solvable by the Cauchy-Kowalevsky theorem, where now
the evolution variable is $y_{2}$ and the initial domain is $\uc^{1}$. 
Notice that equations (C1) will still hold for this extended solution. 
To see this, first differentiate (C2.1) by $y_{1}$. Then substitute, into
the result, equation (C2.2), $k=i$ differentiated by $x_{j}$ and equation
(C2.2), $k=j$ differentiated by $x_{i}$. Thus we obtain 
\[ \dely{2} \left( \dery{\beta_{ij}}{1} - \derx{\alpha_{1j}}{i} + 
\derx{\alpha_{1i}}{j} \right) =0, \]
which tells us that equations (C1) hold for all $y_{2}$ since, 
by the first step, they hold for $y_{2} =0$. Again, the solution will satisfy
also equation (C4.1). This is shown by the same method as in the first step,
except that we use equations (C2.1) instead of (C1). 

The same procedure produces the third and last extension. We have 
ten equations for twelve unknowns. We impose $\alpha_{23}= \alpha_{32}$
and $\alpha_{13}= \alpha_{31}$. Notice that, because of equations (C3.2), 
$k=2$ and (C3.3), $k=1$, we cannot impose $\alpha_{12}= \alpha_{21}$.
So let's treat them as separate unknowns, for the moment. As in the
first and second step we find a solution to the system. Again, we
must show that equations (C1), (C2.1) and (C2.2) are still satisfied.
To prove that (C1) holds we do exactly as in step two
when we proved the same thing, except that we use (C3.1) and (C3.2),
in place of (C2.1) and (C2.2) respectively. We do the same to prove 
that (C2.1) holds, except that we use (C3.1) and (C3.3) and we
differentiate with respect to $y_{2}$ instead of $y_{1}$.
Notice now that from (C1), (C2.1) and (C3.1) we also obtain
(C4.2). To prove that (C2.2) holds, we proceed as follows:
differentiate (C3.2) by $y_{2}$, (C3.3) by $y_{1}$ and (C4.2) by
$x_{k}$. Then, by suitably combining the results, we obtain
 \[\dely{3} \left( \dery{\alpha_{1k}}{2} - \dery{\alpha_{2k}}{1} - \derx{\beta_{12}}{k} \right) = 0, \]
which proves (C2.2). The proof that also (C4.1) holds is just
as in the previous steps. 
It remains to show that $\alpha_{12}= \alpha_{21}$. In fact
it follows from the following:
\[ \begin{array}{lcl}
\dely{3}(\alpha_{12} - \alpha_{21}) &  = &  \dery{\alpha_{32}}{1} + \derx{\beta_{13}}{2} 
    - \dery{\alpha_{31}}{2} - \derx{\beta_{23}}{1} \\
\ & = & - \derx{\beta_{23}}{1} + \derx{\beta_{13}}{2} - \derx{\beta_{12}}{3} \\
\ & = & 0, 
\end{array} \]
where the first equality follows from subtracting (C3.2), $k=2$ and (C3.3), $k=1$;
the second from substituting (C2.2), $k=3$ and using the imposed symmetry
of the other coefficients; the last one is just (C4.1).
The proof is now complete. 
\end{proof}

\begin{observation}{Remark 1.} To prove his more general version of this 
theorem, where the open set $U$ is replaced by any
manifold $M$, Bryant~\cite{bry:calembed} had to use the fact that
every 3-dimensional manifold is parallelizable. His proof then extended
to higher dimensions when $M$ is assumed to be parallelizable. 
To prove Bryant's theorem from our local version, one would need to
understand how to glue solutions obtained from the various coordinate
charts. Accomplishing this might also provide a method to prove the result without
using parallelizability.
\end{observation}

Even though this proof only works locally on a coordinate chart of the given
Riemannian manifold, it is global in the important case of the torus.
\begin{cor} \label{splag:torus}
Let $M$ be the 3-torus with any real-analytic Riemannian metric $g$,
then $(M,g)$ can be isometrically embedded as a special Lagrangian submanifold
of a Calabi-Yau manifold $\mbar$.
\end{cor}
\begin{proof}{Proof.}  We apply Theorem~\ref{splag:embedd} to any standard complexification 
$\uc$ of $M$. We view $g$ as a triply periodic metric tensor in $\reals{3}$, then we
make sure that every choice involved in the three steps of the theorem
is made to be triply periodic in the real part. Solutions will also be triply periodic in the
real part, hence they define a Calabi-Yau structure on $\uc$. Theorem~\ref{bruhat:whitney} 
also ensures that in this way we can describe 
locally all isometric special Lagrangian embeddings of $M$ in some Calabi-Yau 
manifold $\mbar$.
\end{proof}

Given a special Lagrangian torus $M$, one of the questions which arose after
the work of McLean, is whether the family of nearby
special Lagrangian tori, parametrized by the moduli space $\mathcal{M}$, 
actually foliates a neighborhood of $M$ in $\mbar$ (cfr. Introduction). 
This is true in dimension two because
harmonic forms of 2-tori never vanish. In dimension three instead we can construct examples
where this doesn't happen:
\begin{cor}
For any $k \in \nat$, there exist Calabi-Yau manifolds with a special Lagrangian 
3-torus $M$ admitting a harmonic form with $2k$ zeroes, $k$ of which of index $1$ and $k$ of 
index $-1$. Moreover there will be elements of the moduli space of nearby special Lagrangian tori,
arbitrarily close to $M$, intersecting $M$ in at least $2k$ points.
\end{cor}
\begin{proof}{Proof.} In~\cite{calabi:harm} Calabi constructed examples of metrics on the 3-torus which admit
harmonic forms with $k$ zeroes of index $1$ and $k$ of index $-1$. Let $g$ be one of these
metrics and $\theta$ the corresponding harmonic form with zeroes.
As constructed by Calabi, $g$ is not real-analytic, but we can approximate it
(in the $C^{\infty}$ topology) with a real-analytic one $\tilde{g}$. 
The $\tilde{g}$-harmonic form $\tilde{\theta}$ cohomologous to $\theta$ will
also approximate $\theta$ and, by the stability of zeroes of non-zero index,
$\tilde{\theta}$ will have at least the same number of zeroes if the
approximation is precise enough. To the pair $(M, \tilde{g})$ we can then
apply Corollary~\ref{splag:torus} to construct the Calabi-Yau neighborhood $\mbar$.
This proves the first claim.

McLean~\cite{mclean:deform} identified the moduli space of nearby special Lagrangian
tori in $\mbar$ with a three dimensional submanifold $\mathcal{M}$ of $\Gamma(\nu(M))$,
the space of sections of the normal bundle. In fact, given $V \in \mathcal{M}$, the
nearby special Lagrangian torus associated with $V$ is just $M_{V} = \exp_{M}V$.
Via the identification $V \mapsto (JV)^{\flat}$, $\mathcal{M}$ may also be interpreted
as a submanifold of $\Omega^{1}(M)$. As McLean showed, its tangent space at
the zero section is the vector space of harmonic 1-forms. Now let
$\xi(t)$ be a curve in $\mathcal{M}$, viewed in $\Omega^{1}(M)$, such that $\xi(0) = 0$ 
and whose tangent vector at $0$ is $\tilde{\theta}$, the harmonic form
with zeroes.  Then $\lim_{t \rightarrow 0} \xi(t) / t = \tilde{\theta}$ 
in some $C^{k,\alpha}$ topology. Again, by the stability of zeroes of non-zero degree, 
this implies that,  for sufficiently small $t$, 
$\xi(t)$ will have at least the same number of zeroes as $\tilde{\theta}$.
Now if $V(t)$ is the section in  $\Gamma(\nu(M))$ corresponding to $\xi(t)$,
the special Lagrangian submanifold $M_{V(t)}$ will obviously intersect $T$ precisely
at the zeros of $\xi(t)$. This completes the proof.
\end{proof}

\section{Families of special Lagrangian tori} \label{slag:fam}
In the first step of Theorem~\ref{splag:embedd}, 
in the process of finding a solution on $\uc^{1}$, 
we were free to extend arbitrarily almost the entire matrix $A$. This matrix represents the
metric induced by the horizontal slices $U_{t}= \{ y_{1}=t, y_{2}=y_{3}=0 \}$.
So let $A_{t}$ be a choice of this metric for every $t$. We can, for example, ask
the following question: can we choose $A_{t}$ so that every slice $U_{t}$ will
also be special Lagrangian? The following theorem explains when and how
this can be done:
\begin{thm} \label{onepar:constr}
Suppose that $A_{t}$ is a real-analytic one-parameter family of metrics
on $U$.  Then a Calabi-Yau metric can be constructed on $\uc$ so that
each horizontal slice $U_{t}$ is special Lagrangian with metric $A_{t}$ 
if and only if $\det(A_{t})$ does not depend on t and the one form
$(\frac{\partial}{\partial x_{1}})^\flat$ is harmonic w.r.t $A_{t}$ for every $t$.
\end{thm} 
\begin{proof}{Proof.} We use the same notation as in Theorem~\ref{splag:embedd}. 
In particular let the initial metric $g = A_{0}$. In the following,
$x$ will stand short for $(x_{1}, \ldots, x_{3})$
(so, for example, $(x,t,0,0)$ will mean $(x_{1}, \ldots, x_{3},t,0,0)$, in
real coordinates for $\uc$). Imposing the special Lagrangian
condition on the horizontal slices corresponds to
\begin{equation} 
   \left \{ \begin{array}{l}
                \Img \Omega_{(x,t,0,0)}(\delx{1}, \ldots, \delx{3})=0 \\
                B_{t} = 0
              \end{array} 
   \right.
   \label{hor:splag}
\end{equation} 
for all t, where $B_{t}$ is the value of the matrix $B$ on $U_{t}$. A simple 
computation shows that the first one of these holds if and only if:
\[ \Img{\Gamma_{g}}(x,t,0,0) = 0  \]
for all t. Now, since $\Gamma_{g}$ is holomorphic, from this and 
from the Cauchy-Riemann equations we deduce that:
\[ \derx{\Gamma_{g}}{1} (x,t,0,0) = \dery{\Gamma_{g}}{1} (x,t,0,0) = 0, \]
which, by the definition of $\Gamma_{g}$, holds if and only if
\begin{equation} \label{vol:x1}
\derx{\sqrt{g}}{1}(x)=0 
\end{equation}
for all $x \in U$. This is only a condition on the initial data. 
Both conditions in (\ref{hor:splag})
are satisfied if and only if equations (D) and (C1) in the previous section become
\begin{equation} 
    \begin{array}{l}
       \det(A_{t}) = \sqrt{g}(x) \ \text{for all} \ t, \\
       \derx{\alpha_{1j}}{i} - \derx{\alpha_{1i}}{j} =0 \ \ \text{on} \ \uc^{1}.
     \end{array}
     \label{harm:x1}
\end{equation}
It is easy to see that the first equation of (\ref{harm:x1}) together 
with (\ref{vol:x1}) corresponds to the closure of 
$\star(\frac{\partial}{\partial x_{1}})^\flat$ while the second one to the 
closure of $(\frac{\partial}{\partial x_{1}})^\flat$,
so that $(\frac{\partial}{\partial x_{1}})^\flat$ has to be harmonic w.r.t. to $A_{t}$.
The first equation of (\ref{harm:x1}) gives also the independence of $\det(A_{t})$
on $t$. It is also easy to see that these conditions are sufficient to proceed
to the construction of the Calabi-Yau metric on $\uc$ just by following
the second step of Theorem~\ref{splag:embedd}.
\end{proof}

The set of families of metrics $A_{t}$ satisfying the conditions in the Theorem
above is quite rich. In some sense this is a problem because, for example,
one can construct families with metrics degenerating quite badly. On the
other hand we can also easily construct families with behaviors which we expect 
to observe while approaching the singular fibers described
by Gross in \cite{mgross:slfibgeom}. These are are expected to appear in special Lagrangian fibrations 
of compact Calabi-Yau manifolds (cfr. Gross~\cite{mgross:slfibgeom}), but 
some of them have yet to be constructed.

A fairly simple class of such families is the following:
\begin{equation} \label{normal:metrics}
 A_{t}(x_{1}, x_{2}, x_{3}) = \left( \begin{array}{cc}
                                                           e^{u_{t}(x_{1})} & 0 \\
                                                           0 & Q_{t}(x_{1}, x_{2}, x_{3}) 
                                                           \end{array} \right), 
\end{equation}
where $u_{t}$ is any real-analytic function (depending only on $x_{1}$) and $Q_{t}$
is a symmetric, positive definite $2 \times 2$ matrix with real-analytic entries such that
\[ \det(Q_{t}) = e^{-u_{t}(x_{1})} q(x_{2},x_{3}), \]
where $q$ is real-analytic and depending only on $x_{2}$ and $x_{3}$.
If the functions are chosen to be periodic of period $1$ in all three variables, 
$A_{t}$ defines a family of metrics on a three torus, or, if only one
or two are periodic then they are metrics on a cylinder.
The following is the description, in terms of Theorem~\ref{onepar:constr}, of some already 
known examples of one-parameter families of special Lagrangian cylinders: 

\begin{observation}{Example 1.} Suppose that $\sigma: \reals{2} \rightarrow S^{5}$ is
a minimal Legendrian immersion. Then it is known that the cone $C\sigma$ over 
$\sigma(\reals{2})$ is special Lagrangian (cfr. Haskins~\cite{haskins:slcones},
Joyce~\cite{joyce:slspheres}). Also, Haskins and Joyce showed that if 
we consider the one parameter family of curves $\gamma_{t}$ in $\complex{}$
defined by $\gamma_{t} = \{ z \in \complex{} | \Img{z^{3}} = t, \ \arg z \in (0, \pi / 3) \}$
then the one parameter family of manifolds defined by $M_{t} = \gamma_{t} \cdot \sigma(\reals{2})$ 
is smooth, special Lagrangian, asymptotic to the cone $C \sigma$ and degenerating to the cone as $t \rightarrow 0$.
Now parametrize $\gamma_{t}$ by $\gamma_{t}(x_{1}) = (x_{1} + it)^{1/3}$ and
assume, w.l.o.g., that $\sigma$ is conformal. We can thus parametrize each $M_{t}$
by the map $F_{t}: \reals{3} \rightarrow \complex{3}$ given by
\[ F_{t}(x_{1}, x_{2}, x_{3}) = \gamma_{t}(x_{1}) \cdot \sigma(x_2, x_{3}). \]  
It is now easy to see that the metric $A_{t}$ on $M_{t}$, w.r.t. this parametrization, is
\[A_{t} = \left( \begin{array}{ccc}
                 |\dot{\gamma}_{t}|^{2} & 0 & 0 \\
                           0           &  |\gamma_{t}|^{2}f & 0 \\
                           0           &  0            &|\gamma_{t}|^{2}f 
                \end{array}  \right), \]
where $\dot{\gamma}_{t}$ is the derivative w.r.t. to $x_{1}$ and $f ds^{2}$ is the
conformal metric of $\sigma$ (thus $f$ only depends on $x_{2}$ and $x_{3}$). It is
also easy to see that $\det A_{t} = f^{2}/9$, in fact $\dot{\gamma}_{t}\gamma_{t}^{2} =
\frac{1}{3} \frac{d}{dx_{1}}(\gamma_{t}^3) = \frac{1}{3}$. So $A_{t}$ is of the type 
(\ref{normal:metrics}).
One can also check that
\[ \frac{dF_{t}}{dt} = i \cdot dF_{t}(\delx{1}), \]
i.e. that the variational vector field corresponds to the harmonic form $(\delx{1}) ^{\flat}$,
under the identification of the normal bundle with the cotangent bundle. Of course
this is also the case of the families of Theorem~\ref{onepar:constr}. 
As the map $\sigma$ we could for example use the Legendrian, conformal, harmonic maps 
constructed by Haskins~\cite{haskins:slcones} and Joyce~\cite{joyce:slsymmetries}.
\end{observation}

The following two examples show how flexible this construction is. In fact
we choose the family of metrics $A_{t}$,
$t \in [0,t_{1})$, so that the tori start behaving as we would expect if
the family were approaching two
of the singular fibres described by Gross:

\begin{observation}{Example 2.} Choose $Q_{t}$, in (\ref{normal:metrics}), of the following form:
\[ Q_{t} = \left( \begin{array}{cc}
                          1 & 0 \\
                          0 & e^{-u_{t}}
                \end{array}  \right), \]
with $u_{t}$ periodic in $x_{1}$ of period 1. If the following are satisfied:
\[ \begin{array}{c}
           \lim_{t \rightarrow t_{1}} u_{t}(1/2) = +\infty, \\
           \int_{0}^{1} e^{u_{t}(s)/2} ds = 1 \ \text{for all} \ t,
    \end{array} \]
then these metrics describe a family of tori where the 2-cycle $ \{ x_{1} = 1/2 \}$ 
collapses to a circle, while the diameter stays bounded. We expect to
observe a similar behavior near a fibre of type
(2,2) in~\cite{mgross:slfibgeom}.
\end{observation}

\begin{observation}{Example 3.} Now assume
\[ Q_{t} = \left( \begin{array}{cc}
                          e^{v_{t}(x_{1},x_{2})} & 0 \\
                          0                 & e^{-(u_{t}+v_{t})}
                \end{array}  \right). \]
If $u_{t}$ is as in the previous example and $v_{t}$ satisfies:
\[ \begin{array}{c}
           \lim_{t \rightarrow t_{1}} v_{t}(x_{1},1/2) = +\infty \ 
                                                  \text{for all} \ x_{1}, \\
           \int_{0}^{1} e^{v_{t}(x_{1},s)/2} ds = 1 \ \text{for all} \ t
                                                        \ \text{and} \ x_{1},
    \end{array} \]
then also the 2-cycle $\{ x_{2} = 1/2 \}$ will collapse to a circle. 
This is expected to happen while approaching a fibre of type (2,1).
\end{observation}

No example of special Lagrangian fibration containing a fibre of
type $(2,1)$ has been constructed yet.
One approach to the problem of finding one could be to try
to glue this example or similar ones onto a suitable version of the
singular fibre. This though seems, at the moment, a harder problem.
A related question is which of these families can actually
be seen in compact Calabi-Yau's.
We suspect that imposing the curvature of the ambient
manifold to be bounded already provides considerable restrictions
on the types of degenerations occurring in these families.
In fact in Example 3, if we take $v_{t}$ to
depend only on $x_{2}$, one can show that the curvature of
the ambient manifold blows up. For more general choices we
do not know if this still happens. 
We hope to investigate more on these matters in the future.
In the following section we use similar examples to show 
another instance where 3-dimensional special Lagrangian geometry
differs considerably from the 2-dimensional one.

\section{Hitchin's metric is not always Ricci-flat}
Let $\mathcal{M}$ be the local moduli space of the deformations of a
special Lagrangian $n$-torus $M_{0}$ inside an $n$-dimensional
Calabi-Yau manifold $(\mbar, \Omega, \omega)$. 
For each $q \in \mathcal{M}$ denote by $M_{q}$ the special Lagrangian
submanifold corresponding to $q$. As Hitchin~\cite{hitch:msslag} showed, 
$\mathcal{M}$ can be naturally identified with a 
neighborhood of $0$ in $H^1(M_{0}, \reals{})$. 
In the same paper he also proposed the
construction of a Calabi-Yau structure on the so called D-brane
moduli space, i.e. on the manifold
\[ \mathcal{X} = \mathcal{M} \times H^{1}(M_{0}, \reals{}/\integ{}), \]
which according to the SYZ recipe is also a local model for the
Calabi-Yau manifold mirror of $\mbar$. Notice that $\mathcal{X}$
is an $n$-torus fibration over $\mathcal{M}$. Hitchin successfully showed
how to construct naturally an integrable complex structure,
a compatible K\"{a}hler form $\check{\omega}$ and a non vanishing
holomorphic $n$-form $\check{\Omega}$ on $\mathcal{X}$.
This metric is called semi-flat, because it
induces a flat metric on the fibres.
The condition required for these forms to give a Calabi-Yau 
structure is that they are related by the equality
\[ \check{\omega}^n = c \check{\Omega} \wedge \overline{\check{\Omega}} \]
for some constant c. Hitchin proved that this relation 
holds for the proposed forms if
and only if the special Lagrangian submanifolds $M_{q}$ satisfy  a certain
condition. One way to state this condition is the following. Fix a basis
$\Sigma_{1}, \ldots, \Sigma_{n}$ for $H_{1}(M_{0}, \integ{})$. If $\mathcal{M}$ 
is simply connected
then $H_{1}(M_{q}, \reals{})$ can be canonically identified with 
$H_{1}(M_{0}, \integ{})$.
Now, for every $q \in \mathcal{M}$, let $\theta_{1}(q), \ldots, 
\theta_{n}(q)$ be the harmonic 1-forms on $M_{q}$ satisfying
\begin{equation} \label{int:basis}
 \int_{\Sigma_{i}}\theta_{j} = \delta_{ij}.
\end{equation}
Denote by $\inner{\theta_{i}(q)}{\theta_{j}(q)}_{L^{2}}$ the usual 
$L^{2}$ inner product on $\Omega^{1}(M_{q})$ induced by the metric on
$M_{q}$. The condition required then is
that the function
\begin{equation} \label{hitch:cond}
         \Phi :\begin{array}[t]{lcl}
         \mathcal{M} & \rightarrow & \reals{} \\
         q & \mapsto & \det ( \inner{\theta_{i}(q)}{\theta_{j}(q)}_{L^{2}} )
         \end{array}
\end{equation}
is constant on $\mathcal{M}$.

The condition does in fact always hold in the case of special Lagrangian
tori in $K3$ surfaces, see for example Hitchin~\cite{hitch:msclag}. This seemed to give some hope that the same
was true in higher dimensions. Unfortunately it isn't. In this section 
we show that this follows from Theorem \ref{onepar:constr},
which allows us to construct many counterexamples. Had this condition
been true, Hitchin's construction would have provided the first example
of canonical Calabi-Yau structure on the mirror manifold.
 In the final remark we will also show why our counterexamples
fail in dimension 2, as they should. This will highlight what 
goes wrong. So we have:

\begin{cor}
There are 1-parameter families of special Lagrangian tori
along which the function $\Phi$ defined in (\ref{hitch:cond}) 
is not constant.
\end{cor}
\begin{proof}{Proof.} Let $A_{t}$ be a family of metrics on
the standard 3-torus $M =\reals{3}/ \integ{3}$ of the following type:
\[ A_{t} = \left( \begin{array}{lcl}
                     g_{11}(x_{1}, t) & 0 & 0 \\
                     0 & g_{22}(x_{1}, t) & 0 \\
                     0 & 0 & g_{33}(x_{1}, t)
                   \end{array} \right),
\]
with the only condition that $\det(A_{t})= g_{11}g_{22}g_{33} = 1$.
Theorem \ref{onepar:constr} and the comments that follow
show that this family can be realized as a one parameter
family of special Lagrangian submanifolds
of some Calabi-Yau manifold. 
We now show that in general the function $\Phi$ is not constant along 
this family. Choose as basis $\Sigma_{1}, \Sigma_{2}, \Sigma_{3}$ for $H_{1}(M, \integ{})$
the standard one. A computation shows that the forms
\begin{eqnarray*}
\theta_{1} & = & \frac{g_{11}}{\int^{1}_{0}g_{11} dx_{1}} dx_{1}, \\
\theta_{2} & = & dx_{2}, \\
\theta_{3} & = & dx_{3} 
\end{eqnarray*}
are harmonic and they satisfy (\ref{int:basis}) for every $t$. 
Now, since the volume form is just $dx_{1} \wedge dx_{2} \wedge dx_{3}$ 
and the functions given depend only on $x_{1}$ and $t$, we have the
following:
\begin{eqnarray*}
 |\theta_{1}(t)|^{2}_{L^{2}} & = & \frac {1}{ \int_{0}^{1}g_{11}dx_{1}}, \\
 |\theta_{2}(t)|^{2}_{L^{2}} & = & \int_{0}^{1}g^{22}dx_{1}, \\
 |\theta_{3}(t)|^{2}_{L^{2}} & = & \int_{0}^{1}g^{33}dx_{1}, \\
 \inner{\theta_{i}(t)}{\theta_{j}(t)}_{L^{2}} & = & 0
                                        \ \ \text{when} \ i \neq j,
\end{eqnarray*}
where we also used the fact that $g^{ii} = g_{ii}^{-1}$.
Now, using also the condition on the determinant of $A_{t}$, this
implies that
\[ \Phi(t) = \det ( \inner{\theta_{i}(t)}{\theta_{j}(t)}_{L^{2}} ) = 
          \frac{ \int_{0}^{1}g^{22}dx_{1} \int_{0}^{1}g^{33}dx_{1}}
               { \int_{0}^{1}g^{22}g^{33}dx_{1} }, \]
which in general, for arbitrary $g^{22}$ and $g^{33}$ depending
also on $t$, is not constant in $t$.
\end{proof}

\begin{remark}{1} To convince ourselves that these examples show
what goes wrong in dimension 3 and certainly higher, we now demonstrate
why they are not counterexamples in dimension 2, as we expect from
known theory. With slight modifications, one can prove that Theorem 
\ref{onepar:constr} also  holds in dimension 2. Let $A_{t}$ be a
family of metrics on the 2-torus $M = \reals{2}/\integ{2}$
such that $(\delx{1})^\flat$ is harmonic and 
$\det(A_{t})= C(x_{2})$ for every $t$. Then it can be realized
as a one-parameter family of special Lagrangian tori in some
2~-~dimensional Calabi-Yau. We now show that $\Phi$ is constant
along this family. Let $\Sigma_{1}, \Sigma_{2}$ be the standard basis
for $H_{1}(M, \integ{})$. Then, it can be verified that
\begin{eqnarray*}
\theta_{1} & = & \frac{g_{11}(\int^{1}_{0}\sqrt{C} dx_{2}) dx_{1} +
 (g_{12}\int^{1}_{0}\sqrt{C} dx_{2} -\sqrt{C} \int^{1}_{0}g_{12} dx_{2})
dx_{2}}
{\int^{1}_{0}\sqrt{C} dx_{2} \int^{1}_{0}g_{11} dx_{1}}, \\
\theta_{2} & = & \frac{\sqrt{C}}{\int^{1}_{0}\sqrt{C} dx_{2}} dx_{2}
\end{eqnarray*}
are the harmonic $1$-forms satisfying (\ref{int:basis}). Notice that
$\int^{1}_{0}\sqrt{C} dx_{2}$ is just a constant and in fact it
represents the volume of the tori. We can thus assume, w.l.o.g.,
$\int^{1}_{0}\sqrt{C} dx_{2}$=1. 
Also we have that $g^{11} = g_{22}/C$, $g^{22} = g_{11}/C$ and 
$g^{12} = - g_{12}/C$. Using these facts we compute the point-wise
inner product:
\begin{eqnarray*}
|\theta_{1}|^{2}& = &  \frac{g_{11}(g_{11}g_{22} +
                  (g_{12} -\sqrt{C} \int^{1}_{0}g_{12} dx_{2})^{2} - 
                 2g_{12}(g_{12} -\sqrt{C} \int^{1}_{0}g_{12} dx_{2})) }
                  {C(\int^{1}_{0}g_{11} dx_{1})^{2}},                      \\
   & = & \frac{g_{11}(1 + ( \int^{1}_{0}g_{12} dx_{2})^{2})}
                  {(\int^{1}_{0}g_{11} dx_{1})^{2}},                       \\
\inner{\theta_{1}}{\theta_{2}} & = & - \frac{g_{11} \int^{1}_{0}g_{12} dx_{2}}
                                          {\int^{1}_{0}g_{11} dx_{1}} ,    \\
|\theta_{2}|^{2} & = & g_{11}. 
\end{eqnarray*}
Here, to obtain the first equality we have also substituted $g_{11}g_{22} -
g_{12}^{2} = C$. Now, the fact that $( \delx{1})^\flat$ 
is closed implies that $\int^{1}_{0}g_{11} dx_{1}$ and 
$\int^{1}_{0}g_{12} dx_{2}$ are constant. Thus, integrating the
above functions on $M$ yields:
\begin{eqnarray*}
|\theta_{1}|^{2}_{L^{2}} & = & 
           \frac{1 + ( \int^{1}_{0}g_{12} dx_{2})^{2}}
                  {\int^{1}_{0}g_{11} dx_{1}},                       \\
\inner{\theta_{1}}{\theta_{2}}_{L^{2}} & = & -\int^{1}_{0}g_{12} dx_{2},     \\
|\theta_{2}|^{2}_{L^{2}} & = & \int^{1}_{0}g_{11} dx_{1}. 
\end{eqnarray*}
Hence we see that:
\[ \Phi(t) = \det ( \inner{\theta_{i}}{\theta_{j}}_{L^{2}} ) = 1, \]
as we expected. 
\end{remark}

\end{document}